\documentclass{amsart}
\usepackage{amssymb,amsmath,amscd,xy,graphicx,textcomp}
\newtheorem{theorem}{Theorem}[section]
\newtheorem{lemma}[theorem]{Lemma}
\newtheorem{corollary}[theorem]{Corollary}
\newtheorem{definition}[theorem]{Definition}

\newtheorem{remark}[theorem]{\it Remark}

\newtheorem{proposition}[theorem]{Proposition}

\newtheorem{caution}[theorem]{Caution}

\xyoption{arrow}

\xyoption{matrix}

\def\CP{\mathbb{C}P}

\def\C{\mathbb{C}}
\def\R{\mathbb{R}}

\def\br{\mathbf{r}}

\def\tree{\mathcal{T}}

\newcommand{\ZZ}{\mathbb Z}

\begin{document}

\title{Presentations of Semigroup Algebras of Weighted Trees}

\author{Christopher Manon}
\thanks{The author was supported by the NSF FRG grant DMS-0554254.}
\thanks{The author will include these results in his doctoral thesis.}
\date{\today}

\begin{abstract}
We find presentations for subalgebras of
invariants of the coordinate algebras of binary symmetric 
models of phylogenetic trees studied by Buczynska and  Wisniewski
in \cite{BW}. These
algebras arise as toric degenerations of
rings of global sections of weight varieties 
of the Grassmanian of two planes associated to the 
Pl\"ucker embedding, 
and as toric degenerations of rings of invariants of
Cox-Nagata rings.
\end{abstract}

\maketitle

\tableofcontents

\smallskip

\section{Introduction}

Let $\tree$ denote an abstract
trivalent tree with leaves $V(\tree)$,
edges $E(\tree)$, and non-leaf vertices
$I(\tree)$, by trivalent we mean that
the valence of $v$ is three for any $v \in I(\tree)$. 
Let $e_i$ be the unique edge incident to 
the leaf $i \in V(\tree)$.
Let $Y$ be the unique trivalent tree on 
three vertices. For each $v \in I(\tree)$ 
we pick an injective map $i_v: Y \to \tree$,
sending the unique member of $I(Y)$ to $v$. 
We denote the members of $E(Y)$ by $E$, $F$, and $G$.

\begin{definition}
Let $S_{\tree}$ be the graded
semigroup where $S_{\tree}[k]$
is the set of weightings

$$\omega: E(\tree) \to \ZZ_{\geq 0}$$
which satisfy the following conditions.

\begin{enumerate}
\item For all $v \in I(\tree)$ the numbers $i_v^*(\omega)(E)$, 
$i_v^*(\omega)(F)$ and $i_v^*(\omega)(G)$
satisfy 
$$|i_v^*(\omega)(E) - i_v^*(\omega)(F)| \leq i_v^*(\omega)(G) \leq |i_v^*(\omega)(E) + i_v^*(\omega)(F)|$$
These are referred to as the triangle inequalities.  
\item  $i_v^*(\omega)(E) + i_v^*(\omega)(F) + i_v^*(\omega)(G)$ is even. 
\item $\sum _{i \in V} \omega(e_i) = 2k$
\end{enumerate}
\end{definition}

In \cite{SpeyerSturmfels} Speyer and Sturmfels
show that the semigroup algebras
$\C[S_{\tree}]$ may be realized as rings
of global sections for projective embeddings of
the flat toric deformations
of $Gr_2(\C^n)$ where $n = |V(\tree)|,$ for the 
Pl\"ucker embedding. 
This semigroup is also multigraded,
with the grading given by the 
weights $\omega(e_i)$.

\begin{definition}
Let $\br: V(\tree) \to \ZZ_{\geq 0}$
be a vector of nonnegative integers.
Let $S_{\tree}(\br)$ be 
the multigraded subsemigroup 
of $S_{\tree}$ formed by the
pieces  $S_{\tree}[k\br]$.
\end{definition}

It follows from \cite{SpeyerSturmfels}
that graded algebras
$\C[S_{\tree}(\br)]$ are rings of global 
sections for projective embeddings of flat toric deformations of
$Gr_2(\C^n)//_{\br} T$, the weight variety
of the Grassmanian of 2-planes associated
to $\br$, or equivalently $\mathcal{M}_{\br}$, 
the moduli space of $\br$-weighted
points on $\CP^1$.  
In \cite{HowardMillsonSnowdenVakil} this 
degeneration is used to construct presentations of 
the ring of global sections for a projective embedding of $\mathcal{M}_{\br}$, and
it was shown for certain $\br$ and $\tree$ that these algebras
are generated in degree $1$ and have relations generated by 
quadrics and cubics.  This is the starting point for the present paper, 
along with the work of Buczynska and  Wisniewski \cite{BW}, where it was shown
that the algebras associated to the following semigroups
all have the same multigraded Hilbert function, with 
the multigrading defined as it is for $S_{\tree}$. 

\begin{definition}
Let $L$ be a positive integer. 
Let $S_{\tree}^L$ be the graded semigroup
where $S_{\tree}^L[k]$ is the set of weightings
$\omega$
of $\tree$ which satisfies the same conditions
as $S_{\tree}[k]$ with the addition assumption
that
$$i_v^*(\omega)(E) + i_v^*(\omega)(F) + i_v^*(\omega)(G) \leq 2kL.$$
This is referred to as the level condition.
\end{definition}

In \cite{SturmfelsXu}, 
Sturmfels and Xu construct the multigraded
Cox-Nagata ring $R^G(L),$ which functions as an analogue
of $Gr_2(\C^n)$ in that it can be flatly deformed
to each $\C[S^L_{\tree}]$. The analogue of the 
weight varieties in this context are the multigrade $\br$
Veronese subrings of $R^G(L)$, denoted $R^G(L)_{\br}$. 

\begin{definition}
Let $L$ be a positive integer. 
Let $S_{\tree}^L(\br)$ be the multigraded
subsemigroup of $S_{\tree}^L$
of summands with multigrade $k\br$. 
\end{definition} 

\begin{remark}
The multigraded Hilbert functions
of $\C[S_{\tree}^L]$ and $\C[S_{\tree}^L(\br)]$ are closely
related to the Verlinde Formula for $SU(2)$ 
(see \cite{BW} and \cite{SturmfelsXu}). Indeed,
the reader may notice that 
the defining conditions for $S_{\tree}^L$
and $S_{\tree}^L(\br)$
are given by Quantum Clebsch-Gordon
Rules for $SU(2)$, whereas the defining conditions
for $S_{\tree}$ and $S_{\tree}(\br)$ are
classical $SU(2)$ Clebsch-Gordon Rules.
\end{remark}

It follows from results in \cite{SturmfelsXu} that
$\C[S_{\tree}^L(\br)]$ is a toric
deformation of $R^G(L)_{\br}$. 
In this paper we construct presentations
for a large class of the rings $\C[S_{\tree}^L(\br)]$.
The techniques used are such that the same results
immediately hold for $\C[S_{\tree}(\br)]$ as well, in 
particular we give a different proof of a fundamental
result of \cite{HowardMillsonSnowdenVakil} on a presentation
of these rings.

\subsection{Statement of Results}

We now state the main results of the
paper. When two leaves are both connected
to a common vertex, we say they are 
paired to each other.  A leaf that has 
no pair is called a lone leaf. 

\begin{definition}
We call the triple $(\tree, \br, L)$
admissible if $L$ is even, $\br(i)$ is even
for every lone leaf $i$, and 
$\br(j) + \br(k)$ is even 
for all paired leaves $j$, $k$.
\end{definition}

\begin{remark}
Admissability is actually
not very restrictive. Note that the assumption
that $\br$ has an even total sum gurantees that 
we may find a $\tree$ such that $(\tree, \br, L)$
is admissible for any even $L$.  This is important
for constructing presentations of 
$R^G(L)^T_{\br}$, since this ring always has a flat
deformation to $\C[S_{\tree}^L(\br)]$ for some 
admissible $(\tree, \br, L)$. 
Also note that the second Veronese subring
of $\C[S_{\tree}^L(\br)]$ is the semigroup
algebra associated to $(\tree, 2\br, 2L)$, 
which is always admissible.  
\end{remark}

\begin{theorem}\label{gen}
For $(\tree, \br, L)$ admissible with $L > 2$,
$\C[S_{\tree}^L(\br)]$ is generated in degree $1$.
\end{theorem}

\begin{theorem}\label{rel}
For $(\tree, \br, L)$ admissible with $L > 2$, 
$\C[S_{\tree}^L(\br)]$ has relations
generated in degree at most $3$. 
\end{theorem}
As a corollary we get the same results for $S_{\tree}(\br)$
when $(\tree, \br)$ satisfy admissibility conditions.
These theorems will be proved in 
sections $2$, $3$, and $4$.  In section $5$
we will look at some special cases, and investigate
what can go wrong when $(\tree, \br, L)$ is not
an admissible triple.  

\subsection{Outline of techniques}

To prove Theorems \ref{gen} and \ref{rel} we
use two main ideas.  First, we employ the following
trivial but useful observation. 

\begin{proposition}
Let $(\tree, \br, L)$ be admissible, 
then for any weighting $\omega \in S_{\tree}^L(\br)$, 
$\omega(e)$ is an even number when $e$ is not
an edge connected to a paired leaf. 
\end{proposition}
This allows us to drop the parity
condition that 
$i_v^*(\omega)(E) + i_v^*(\omega)(F) + i_v^*(\omega)(G)$
is even by forgetting the paired leaves
and halving all remaining weights.  

\begin{definition}
Let $c(\tree)$ be the subtree
of $\tree$ given by forgetting
all edges incident to paired leaves.
\end{definition} 

\begin{figure}[htbp]
\centering
\includegraphics[scale = 0.2]{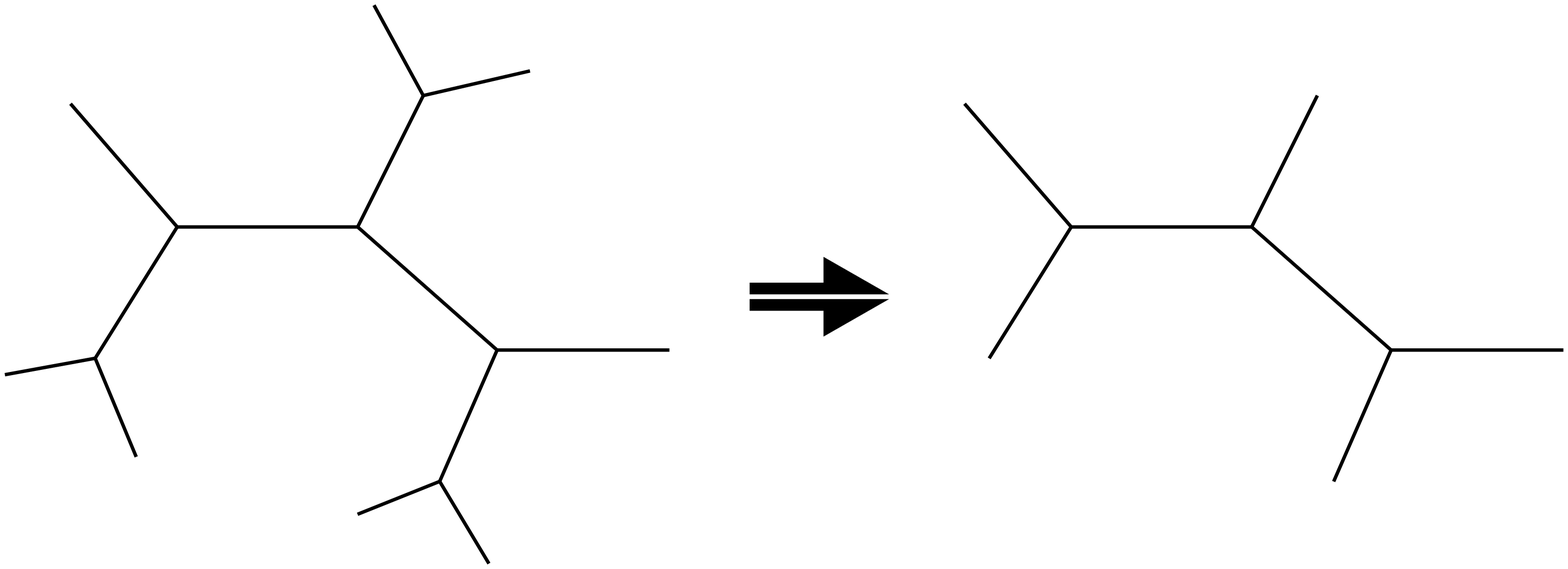}
\caption{Clipping the paired leaves of $\tree$}
\label{fig:Clipping}
\end{figure}

\begin{definition}
Let $U_{c(\tree)}^L(\br)$
be the graded semigroup of weightings on $c(\tree)$
such that the members of $U_{c(\tree)}^L(\br)[k]$
satisfy the triangle inequalities,
the new level condition $i_v^*(\omega)(E) + i_v^*(\omega)(F) + i_v^*(\omega)(G) \leq L$, 
and the following conditions. 
\begin{enumerate}
\item $\omega(e_i) = k \frac{\br(i)}{2}$ for $i$ a lone leaf of $\tree$.
\item $\frac{k|\br(j) - \br(k)|}{2} \leq \omega(e) \leq \frac{k|\br(j) + \br(k)|}{2}$ 
for $e$ the unique edge of $\tree$ connected to the vertex
which is connected to the paired leaves $i$ and $j$.
\end{enumerate}
\end{definition}
Also, let $U_{c(\tree)}^L$ be the graded semigroup 
of weightings which satisfy the triangle inequalities
and the new level condition for $L$. The following is a 
consequence of these definitions. 

\begin{proposition}\label{sum2gone}
For $(\tree, \br, L)$ admissible, 

$$U_{c(\tree)}^L(\br) \cong S_{\tree}^L(\br)$$
as graded semigroups. 
\end{proposition}
We refer to the graded semigroup of weightings
on $Y$ which satisfy the triangle inequalities
and the new $L$ level condition as $U_Y^L$. 
The next main idea is to undertake the analysis 
of $U_{c(\tree)}^L(\br)$ by first considering 
the weightings $i_v^*(\omega) \in U_Y^L$.  
After constructing a pertinant
object in $U_Y^L$, like a factorization or relation, 
we ``glue'' these objects
back together along edges shared  
by the various $i_v(Y)$ with what amounts
to a fibered product of graded semigroups.  
This is remniscent of the theory of
moduli of orientable surfaces, 
where structures on a surface of high
genus can be glued together
from structures on three-punctured spheres over
a pair-of-pants decomposition.  The
reason for this resemblance is not entirely
accidental, see \cite{HowardManonMillson}. 
We obtain information about $U_Y^L$ 
by studying the following polytope. 

\begin{remark}
In \cite{BW}, Buczynska and  Wisniewski used
more or less the same idea. 
They prove facts about $S_{\tree}^L$ by viewing it
as a fibered product of copies of $S_Y^L$. 
\end{remark}

\begin{definition}
Let $P_3(L)$ be the convex hull
of $(0, 0, 0)$, $(\frac{L}2, \frac{L}2, 0)$, $(\frac{L}2, 0, \frac{L}2)$, 
and $(0, \frac{L}2, \frac{L}2)$. 
\end{definition} 

\begin{figure}[htbp]
\centering
\includegraphics[scale = 0.3]{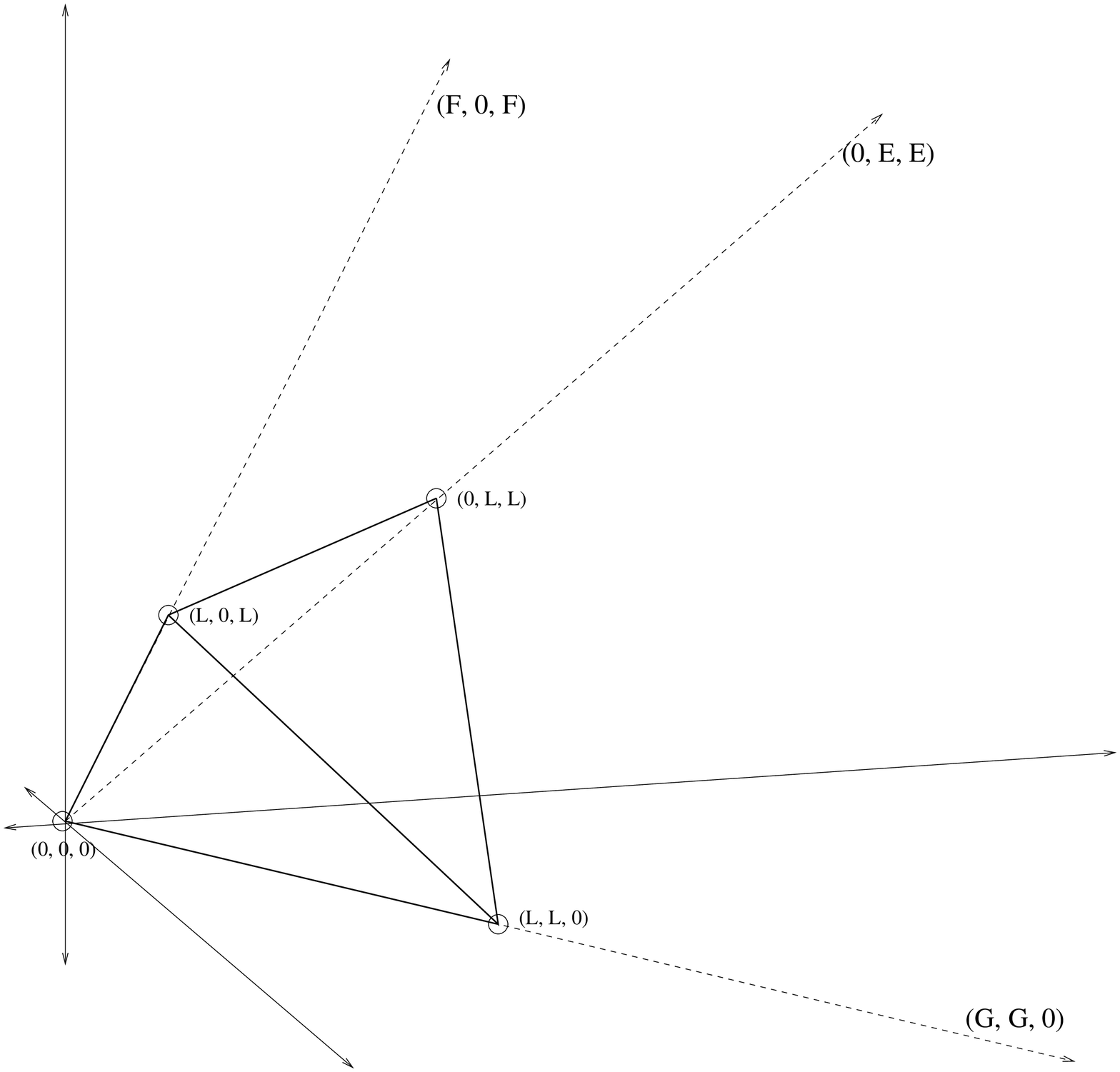}
\caption{$P_3(2L)$}
\label{fig:Polytope}
\end{figure}
\begin{figure}[htbp]
\centering
\includegraphics[scale = 0.3]{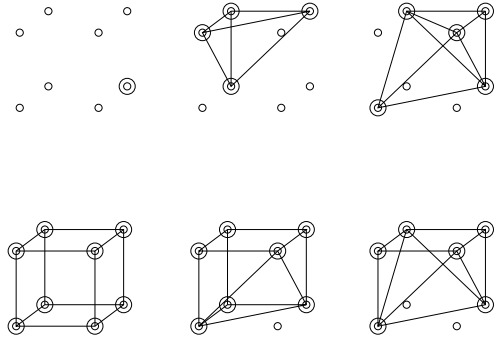}
\caption{Cube Polytopes}
\label{fig:cubes}
\end{figure}

The graded semigroups of lattice points
for $P_3(L)$ is $U_Y^L$. 
By a lattice equivalence of polytopes $P$, $Q \subset \R^n$
with respect to a lattice $\Lambda \subset \R^n$ we mean a composition
of translations by members of $\Lambda$ and members of $GL(\Lambda) \subset GL_n(\R)$ 
which takes $P$ to $Q$.  
If $P$ and $Q$ are lattice
equivalent it is easy to show that they have isomorphic
graded semigroups of lattice points. 
When $L$ is an even integer (admissibility condition) the interesection
of this polytope with any translate of the unit cube in $\R^3$, 
is, up to lattice equivalence, one of the polytopes shown in figure \ref{fig:cubes}. 

Each of these polytopes is normal, 
and the relations of their
associated semigroups are generated
in degree at most $3$.  In sections $3$
and $4$ we will lift these
properties to $U_{c(\tree)}^L(\br)$, 
and therefore $S_{\tree}^L(\br)$ for $(\tree, \br, L)$
admissible.  Facts about the six polytopes
above also allow us to carry out a more
detailed investigation into the properties
of the semigroups $S_{\tree}^L(\br)$ in section $5$, for
example they allow us to show the
redundancy of the cubic
relations for certain $(\tree, \br, L)$.\\

I would like to thank John Millson for introducing me to this 
problem, Ben Howard for many useful and encouraging conversations
and for first introducing me to the commutative algebra of
semigroup rings, Larry O'Neil for several useful conversations
on the cone of triples which satisfy the triangle inequalities, 
and Bernd Sturmfels for introducing me to Graver bases and shortening the proof
of Theorems and \ref{cubegen} and \ref{cuberel}.

\section{The Cube Semigroups}

In this section we will prove that
the intersection of any translate of the unit cube of $\R^3$
with $P_3(2L)$ produces
a normal polytope whose semigroup of lattice
points has relations generated in degree
at most $3$ when $L$ is an integer. 
Let $P_3$ be the cone of triples of
nonnegative integers which satisfy
the triangle inequalities, and let 

$$C(m_1, m_2, m_3) = conv\{ (m_1 + \epsilon_1, m_2 + \epsilon_2, m_3 + \epsilon_3)| \epsilon_i \in \{0, 1\} \}$$
We want classify the polytopes which have the
presentation $C(m_1, m_2, m_3)\cap P_3$, since $P_3$ is symmetric 
we may assume that $(m_1, m_2, m_3)$ is ordered by magnitude with $m_3$ the largest.
We will keep track of the triangle inequalities with the quantities $n_i = m_j + m_k - m_i$. For
a point $(m_1, m_2, m_3)$ to be in $P_3$ is equivalent to $n_i \geq 0$ for each $i$. 
Immediately we have the following inequalities.
$$n_1 \geq n_2 \geq n_3, n_2 \geq 0$$
If $n_3 < -2$ then no member of $C(m_1, m_2, m_3)$ can belong
to $P_3$.  If $n_3 \geq -2$ then there are six distinct possibilities, we list each case
along with the members of $C(m_1, m_2, m_3)\cap P_3 - (m_1, m_2, m_3)$.
\\

\begin{tabular}[htbp]{|r|r|}
\hline
$Condition$&$C(m_1, m_2, m_3) \cap P_3 - (m_1, m_2, m_3)$\\
\hline
$n_3 = -2$&$(1, 1, 0)$\\
$n_3 = -1$&$(1, 1, 0), (0, 1, 0), (1, 0, 0), (1, 1, 1)$\\
$n_1 = n_2 = n_3 = 0$&$ (1, 1, 0), (0, 1, 1), (1, 0, 1), (1, 1, 1), (0, 0, 0)$\\
$n_1 > 0, n_2 = n_3 = 0$&$(1, 1, 0), (0, 1, 1), (1, 0, 1), (1, 1, 1), (0, 0, 0), (0, 0, 1)$\\
$n_1, n_2 > 0, n_3 = 0$&$(1, 1, 0), (0, 1, 1), (1, 0, 1), (1, 1, 1), (0, 0, 0), (0, 0, 1), (0, 1, 0)$\\
$n_i > 0$& all points \\
\hline
\end{tabular}
\\

The figure below illustrates these arrangements.

\begin{figure}[htbp]
\centering
\includegraphics[scale = 0.3]{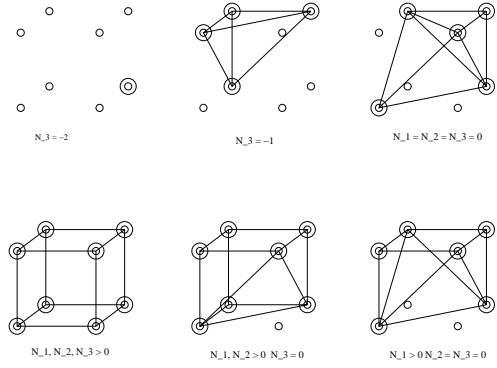}
\caption{ Primitive cube semigroups}
\label{fig:Cube}
\end{figure}

\begin{figure}[htbp]
\centering
\includegraphics[scale = 0.3]{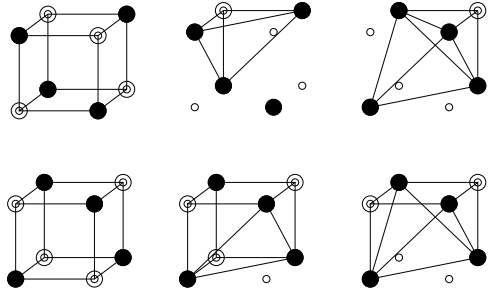}
\caption{Cube semigroups with the lattice $v_1 + v_2 +v_3 = 0$ $mod$ $2$}
\label{fig:Parity}
\end{figure}

Now we will see what happens when we
intersect  $P_3$ with the half plane
$m_1 + m_2 + m_3 \leq 2L$ to get $P_3(2L)$.
The reader may want to refer to figure 
\ref{fig:Parity} for this part. 
The convex set $C(m_1, m_2, m_3) \cap P_3(2L)$ can be
one of the above polytopes (up to lattice equivalence), 
or one of them intersected with the half plane
$m_1 + m_2 + m_3 \leq 2L$. 
Note that a vertex $v$ in $C(m_1, m_2, m_3) \cap P_3(2L)$
lying on a facet of $P_3$ necessarily
satisfies $v_1 + v_2 +v_3 = 0$ $mod$ $2$. 
In Figure \ref{fig:Parity} these points are
colored black.

\begin{figure}[htbp]
\centering
\includegraphics[scale = 0.3]{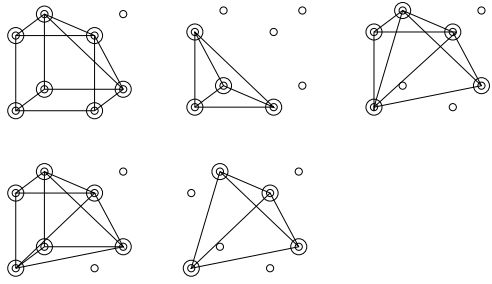}
\caption{New Possibilities for  $C(m_1, m_2, m_3) \cap P_3(2L)$ }
\label{fig:LevelPoly}
\end{figure}

The hyperplane $m_1 + m_2 + m_3 = 2L$ must intersect
these polytopes at collections of three
black points. If we assume that the lower left
corner is $(0, 0, 0)$, these points have coordinates
$\{(1, 1, 0)$, $(1, 0, 1), (0, 1, 1)\}$, or 
$\{(1, 0, 0), (0, 1, 0), (0, 0, 1)\}$. 
Figure \ref{fig:Parity} represents the 
new possibilities for 
$C(m_1, m_1, m_3) \cap P_3(2L) - (m_1, m_2, m_3)$. 
The polytope pictured lower
center in Figure \ref{fig:LevelPoly}
is the only case which is not lattice equivalent 
to one pictured in Figure \ref{fig:Cube}.
It is rooted at $(0, 0, 0)$ and occurs only
when $L = 1$ (level condition is $2$). The point
$(1, 1, 1)$ in its second
Minkowski sum cannot be expressed as the sum of two 
lattice points of degree one, so this is not a normal 
polytope.  This is the reason we
stipulate that $L > 2$ in Theorem \ref{gen}. 
Now we analyze each $C(m_1, m_2, m_3)\cap P_3(2L)$.
Since lattice equivalent polytopes
have isomorphic semigroups of lattice points, 
it suffices to investigate 
the polytopes listed in Figure \ref{fig:Cube}.

\begin{caution}
In [BW], Buczynska and  Wisniewski study
a normal polytope with the same vertices as the non-normal
polytope mentioned above. 
This is possible because they are using the the
lattice $v_1 + v_2 + v_3 = 0$ $mod$ $2$, not
the standard lattice.
\end{caution}

We make use of the computational algebra package
4ti2, \cite{4ti2} to compute the $Graver$ $basis$ 
of the toric ideal of the unit 3-cube.

$$
\begin{tabular}[htbp]{rr}
(1,0,0) + (1,1,1) = (1,0,1) + (1,1,0)& (0,1,0) + (1,1,1) = (0,1,1) + (1,1,0)\\
(0,0,0) + (1,1,1) = (0,0,1) + (1,1,0)& (0,0,1) + (1,1,1) = (0,0,1) + (1,0,1)\\
(0,0,0) + (1,1,1) = (0,1,0) + (1,0,1)& (0,0,1) + (1,1,0) = (0,1,0) + (1,0,1)\\
(0,0,0) + (1,1,1) = (0,1,1) + (1,0,0)& (0,0,1) + (1,1,0) = (0,1,1) + (1,0,0)\\
(0,1,0) + (1,0,1) = (0,1,1) + (1,0,0)& (0,0,0) + (1,1,0) = (0,1,0) + (1,0,0)\\
(0,0,0) + (1,0,1) = (0,0,1) + (1,0,0)& (0,0,0) + (0,1,1) = (0,0,1) + (0,1,0)\\
\end{tabular}
$$

$$
\begin{tabular}[htbp]{r}
(0,1,0) + (1,0,0) + (1,1,1) = (0,0,1) + (1,1,0) + (1,1,0)\\
(0,0,0) + (1,1,1) + (1,1,1) = (0,1,1) + (1,0,1) + (1,1,0)\\
(0,0,1) + (1,0,0) + (1,1,1) = (0,1,0) + (1,0,1) + (1,0,1)\\
(0,0,1) + (0,0,1) + (1,1,0) = (0,0,0) + (0,1,1) + (1,0,1)\\
(0,0,0) + (0,1,1) + (1,1,0) = (0,1,0) + (0,1,0) + (1,0,1)\\
(0,0,0) + (1,0,1) + (1,0,1) = (0,1,1) + (1,0,0) + (1,0,0)\\
(0,0,1) + (0,1,0) + (1,1,1) = (0,1,1) + (0,1,1) + (1,0,0)\\
(0,0,0) + (0,0,0) + (1,1,1) = (1,0,0) + (0,1,0) + (0,0,1)\\
\end{tabular}
$$

Operating on this set of monomials, one can
show that the toric ideal of every sub-polytope 
of the unit 3-cube which is not a simplex has a square-free
Gr\"obner basis.  This, combined with the fact
that the sub-polytopes with $n_3 = -2$ and $-1$
are unimodular simplices shows the following theorem. 

\begin{theorem}\label{cubegen}
Let $L \neq 1$, then for all $(m_1, m_2, m_3)$, 
if $C(m_1, m_2, m_3)\cap P_3(2L)$
is non-empty, then it is normal. 
\end{theorem}

\begin{remark}\label{flooredge}
This theorem implies, among
other things, that if 
$\omega \in U_Y(2L)[k]$, then

$$\omega = \sum_{i =1}^k W_i$$
for $W_i \in P_3(2L)$
with the property that each

$$W_i = X + (\epsilon_1, \epsilon_2, \epsilon_3)$$
with $\epsilon_j \in \{0, 1\}$ for all $i$ for a
fixed $X \in \R^3$. 
It is easy to show that 

$$X = (\lfloor \frac{\omega(E)}{k} \rfloor, \lfloor \frac{\omega(F)}{k} \rfloor, \lfloor \frac{\omega(G)}{k} \rfloor)$$
Therefore each $W_i$ is $(\frac{\omega(E)}k, \frac{\omega(F)}k, \frac{\omega(G)}k)$
with either floor or ceiling applied to each entry. 
\end{remark}

Now we move on to relations, Let $S(m_1, m_2, m_3)$
be the semigroup of lattice points for 
$C(m_1, m_2, m_3) \cap P_3(2L) - (m_1, m_2, m_3),$
once again it suffices to treat the cases
represented in Figure \ref{fig:Cube}.
 
\begin{theorem}\label{cuberel}
All relations for the semigroup
$S(m_1, m_2, m_3)$ are reducible
to quadrics and cubics. 
\end{theorem}

\begin{proof}
By Proposition 4.13 of \cite{Sturmfels}, a Graver basis
for any subpolytope $P$ of the unit 3-cube is obtained by 
taking the members of the Graver basis of the unit
3-cube which have entries in the lattice points of $P$.
Since these are all quadrics and cubics, we are done. 
\end{proof}

Up to equivalence and after accounting
for redundancy, all relations are of the form
$$(1, 0, 0) + (0, 1, 0) = (1, 1, 0) + (0, 0, 0)$$
$$(1, 0, 1) + (0, 1, 0) = (1, 1, 1) + (0, 0, 0)$$
$$(1, 0, 1) + (1, 1, 0) = (1, 1, 1) + (1, 0, 0)$$
$$(1, 1, 1) + (1, 1, 1) + (0, 0, 0) = (1, 1, 0) + (1, 0, 1) + (0, 1, 1),$$
with the last one the only degree $3$ relation, we refer to it as the 
``degenerated Segre Cubic'' (see \cite{HowardMillsonSnowdenVakil}).

\section{Proof of Theorem \ref{gen}}
In this section we use Theorem
\ref{cubegen} to prove that 
$U_{c(\tree)}^L(\br)$ is generated
in degree $1$, which then proves 
Theorem \ref{gen}. 
For each $v \in I(\tree)$
we have the morphism of graded semigroups

$$i_v^*: U_{c(\tree)}^L(\br) \to U_Y^L.$$
Given a weight $\omega \in U_{c(\tree)}^L(\br)$
we factor $i_v^*(\omega)$ for each $Y \subset c(\tree)$
using Theorem \ref{cubegen}. Then, special properties
of the weightings obtained by this procedure
will allow us to glue the factors of the 
$i_v^*(\omega)$ back together along common edges
to obtain a factorization
of $\omega$. First we must make sure
that our factorization procedure
does not disrupt the conditions at the
edges of $c(\tree)$.

\begin{lemma}
Let $\omega \in U_{c(\tree)}^L(\br)[k]$, 
and let $v \in I(\tree)$ be connected
to a leaf of $c(\tree)$, at $E$. Then
if $i_v^*(\omega) = \eta_1 + \ldots +\eta_k$
is any factorization of $i_v^*(\omega)$
with $\eta_i \in C(\lfloor \frac{i_v^*(\omega)(E)}k \rfloor, \lfloor \frac{i_v^*(\omega)(F)}k \rfloor, \lfloor \frac{i_v^*(\omega)(G)}k \rfloor)$
Then $\eta_i(E)$ satisfies the appropriate edge
condition for elements in  $U_{c(\tree)}^L(\br)[1]$.
\end{lemma} 

\begin{proof}
If $E$ is attached to a lone
leaf of $\tree$ then $i_v^*(\omega)(E) = k \br(e)$ for 
$i_v(E) = e$,  $e \in V(\tree)$.  By Remark \ref{flooredge}

$$\eta_i(E) = \lfloor \br(e) \rfloor = \br(e)$$
or

$$\eta_i(E) = \lfloor \br(e) \rfloor + 1 = \br(e) + 1$$
Since $\sum_{i = 1}^k \eta_i(E) = k \br(e)$ 
we must have $\eta_i(E) = \br(e)$ for all $i$. 
If $E$ is a stalk of paired leaves $i$ and $j$ in $\tree$
then we must have

$$k\frac{|\br(i) - \br(j)|}2 \leq \omega_Y(E) \leq k\frac{|\br(i) + \br(j)|}2$$
Note that both bounds are divisible by $k$. 
Since floor preserves lower bounds we have

$$\frac{|\br(i) - \br(j)|}2 \leq \lfloor \frac{i_v^*(\omega)(E)}k \rfloor,$$
and since ceiling preserves upper bounds
we have

$$\lceil \frac{i_v^*(\omega)(E)}k \rceil \leq \frac{|\br(i) + \br(j)|}2.$$
Therefore each $\eta_i$ satisfies

$$\frac{|\br(i) - \br(j)|}2 \leq \eta_i(E) \leq \frac{|\br(i) + \br(j)|}2$$
\end{proof}
Now that we can safely
use Theorem \ref{cubegen} with
each $i_v^*:U_{c(\tree)}^L(\br) \to U_Y^L$, we can 
see about gluing these factors
together along common edges. 

\begin{definition}
We say a set of nonnegative
integers $\{X_1, \ldots, X_n\}$
is balanced if $| X_i - X_j| = 1$
or $0$ for all $i$, $j$. 
\end{definition}
The following is a very useful 
lemma, its proof is left to the reader.  

\begin{lemma}\label{bsets}
If two sets $\{X_1, \ldots, X_n\}$
and $\{Y_1, \ldots, Y_m\}$ are balanced, 
have the same total sum, 
and $n = m$, then they are the same set.
\end{lemma}

\begin{proposition}
The semigroup $U_{c(\tree)}^L(\br)$ is 
generated in degree $1$.
\end{proposition}

\begin{proof}
Recall that by Remark \ref{flooredge}, 
for any edge $E \in Y$ 
the edge weights of a factorization
$i_v^*(\omega) = \eta_1 + \ldots \eta_k$
satisfy $\eta_i(E) = \lfloor \frac{i_v^*(\omega)(E)}k \rfloor$
or $\lceil \frac{i_v^*(\omega)(E)}k \rceil$.
Take any two $v_1$, $v_2$ which share
a common edge $E$ in $c(\tree)$. Let $\omega \in U_{c(\tree)}^L(\br)[k]$
and let $\{\eta^1_1, \ldots, \eta^1_k\}$
and $\{\eta^2_1, \ldots, \eta^2_k\}$ be
factorizations of $i_{v_1}^*(\omega)$ and $i_{v_1}^*(\omega)$
respectively. Then the sets
$\{ \eta^1_1(E), \ldots, \eta^1_k(E)\}$
and $\{\eta^2_1(E), \ldots, \eta^2_k(E)\}$
are balanced and have the same sum, 
so by Lemma \ref{bsets} they are the same set. 
We may glue factors $\eta_i^1$
and $\eta_j^2$ when $\eta_i^1(E) = \eta_j^2(E)$, 
the above observation gurantees that any $\eta_i^1$
has an available partner $\eta_j^2$. The proposition
now follows by induction on the number of
$v \in I(c(\tree))$. This implies Theorem \ref{gen}. 
\end{proof}

\section{Proof of Theorem \ref{rel}}

In this section we show how
to get all relations in $U_{c(\tree)}^L(\br)$
from those lifted from $U_Y^L$.
The procedure follows the same pattern
as the proof of Theorem \ref{gen}. We consider
the image of a relation 
$\omega_1 + \ldots + \omega_n = \eta_1 + \ldots + \eta_n$
under a map $i_v^*: U_{c(\tree)}^L(\br) \to U_Y^L$, using
Theorem \ref{cuberel} we convert this to a trivial relation 
using relations of degree at most $3$. We then give a recipe 
for lifting each of these relalations 
back to $U_{c(\tree)}^L(\br)$. The result
is a way to convert $\omega_1 + \ldots + \omega_n = \eta_1 + \ldots + \eta_n$
to a relation which is trivial over the trinode $v$ using
quadrics and cubics. In this way we take a general relation to a 
trivial relation one $v \in I(c(\tree))$ at a time.

\begin{definition}
A set of degree $1$ elements 
$\{\omega_1, \ldots, \omega_k\}$ in $U_{c(\tree)}^L(\br)$ is
called Balanced when the set 
$\{\omega_1(E), \ldots, \omega_k(E)\}$
is balanced for all $E \in c(\tree)$. 
A relation $\omega_1 + \ldots + \omega_k = \eta_1 + \ldots + \eta_k$
in $U_{c(\tree)}^L(\br)$ is called Balanced when 
$\{\omega_1, \ldots, \omega_k\}$ and $\{\eta_1, \ldots, \eta_k\}$
are balanced. 
\end{definition}
The following lemmas say that 
we need only consider balanced relations. 

\begin{lemma}\label{convert}
Any set of nonegative integers
$S = \{X_1, \ldots, X_n\}$ can 
be converted to a balanced set $T = \{Y_1, \ldots, Y_n\}$ with $\sum_{i = 1}^n Y_i = \sum_{i = 1}^n X_i$ 
by replacing a pair $X_i$ and $X_j$ with $\lfloor \frac{X_i + X_j}2 \rfloor$
and $\lceil \frac{X_i + X_j}2 \rceil$ a finite number of times.
\end{lemma}

\begin{proof}
Let $d(S)$ be the difference between the
maximum and mininum elements of $S$.  It is 
clear that with a finite number of exchanges 
$$\{X_i, X_j\} \to \{\lfloor \frac{X_i + X_j}2 \rfloor, \lceil \frac{X_i + X_j}2 \rceil\}$$
We get a new set $S'$ with $d(S) > d(S')$, unless
$d(S) = 1$ or $0$.  Since this happens if and only of $S$
is balanced, the lemma follows by induction.
\end{proof}

\begin{lemma}
Let 
$$\omega_1 + \ldots + \omega_k = \eta_1 + \ldots + \eta_k$$
be a relation in $U_{c(\tree)}^L(\br)$ then it can 
be converted to a balanced relation 
$$\omega_1' + \ldots + \omega_k' = \eta_1' + \ldots + \eta_k'$$
using only degree $2$ relations. 
\end{lemma}

\begin{proof}
First we note that using
the proof of Theorem \ref{gen}
we can factor the weighting $\omega_1 + \omega_2$
into $\omega_1' + \omega_2'$ so that $\{\omega_1', \omega_2'\}$
is balanced.  Using this and Lemma \ref{convert} we can find
$$\omega_1' + \ldots + \omega_k' = \omega_1 + \ldots + \omega_k$$
such that the set $\{\omega_1'(E), \ldots, \omega_k'(E)\}$  
is balanced for some specific $E$, using only degree $2$ relations. 
Observe that if $\{\omega_1(F), \ldots, \omega_k(F)\}$ is balanced
for some $F$, the same is true for  
$\{\omega_1'(F), \ldots, \omega_k'(F)\}$, 
after a series of degree $2$ applications of \ref{gen}
as above.  This shows that we may inductively
convert $\{\omega_1, \ldots, \omega_k\}$
to $\{\omega_1', \ldots, \omega_k'\}$ with the
property that $\{\omega_1'(E), \ldots, \omega_k'(E)\}$
is a balanced set for all edges $E$, using only
degree $2$ relations.  Applying the same procedure
to the $\eta_i$ then proves the lemma.  
\end{proof}
The next lemma shows how we lift
a balanced relation in $U_Y^L$ to one
in $U_{c(\tree)}^L(\br)$. 

\begin{lemma}\label{lifting}
Let $\{\omega_1 \ldots \omega_k\}$ be a balanced set of 
elements in $U^{L}_{c(\tree)}(\br)$.  Let
$i_v^*(\omega_1) + \ldots + i_v^*(\omega_k) = \eta_1 + \ldots + \eta_k$ 
be a degree $k$ relation the appropriate $S(m_1, m_2, m_3) \subset U_Y^L$. 
Then the $\eta_i$ may be lifted to weightings of $c(\tree)$ giving 
a relation of degree $k$ in $U_{c(\tree)}^L(\br)$ which agrees with
the relation above when $i_v^*$ is applied, and is a permutation of 
$i_{v'}^*(\omega_1) \ldots i_{v'}^*(\omega_N)$ for $v' \neq v$.
\end{lemma}

\begin{figure}[htbp]
\centering
\includegraphics[scale = 0.4]{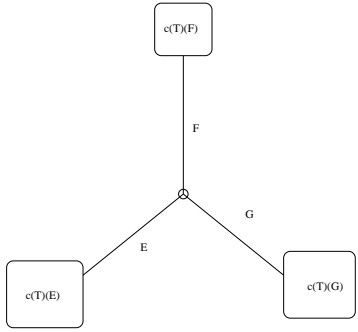}
\caption{ Component subtrees about a vertex}
\label{fig:Components}
\end{figure}

\begin{proof}
Let $c(\tree)(E)$ be the unique connected subtrivalent 
tree of $c(\tree)$ which includes $v$ and has the property that any
path $\gamma \subset c(\tree)(E)$ with endpoints at a vertex 
$v' \neq v$ in $c(\tree)(E)$ and $v$ includes the edge $E$
(see Figure \ref{fig:Components}),
define $c(\tree)(F)$ and $c(\tree)(G)$ in the same way. 
To make $\eta_1' \ldots \eta_k'$ over $c(\tree)$, note 
that the set $\{i_{c(\tree)(E)}^*(\omega_i)(E)\}$ is the same as 
the set $\{\eta_i(E)\}$, because they are both balanced
sets with the same sum and the same number of elements, so we 
may glue these weightings together to make a tuple over $c(\tree)$. 
\end{proof}

If we are given a relation
$$\omega_1 + \ldots + \omega_k = \eta_1 + \ldots + \eta_k$$
with both sides balanced, we may use relations
in the appropriate $S(m_1, m_2, m_3)$ to convert
$\{\omega_1, \ldots, \omega_k\}$ to 
$\{\eta_1, \ldots, \eta_k\}$ one $v \in I(c(\tree))$ at a time. 
This leads us to the following proposition. 

\begin{proposition}
Let $N$ be the maximum degree of relations
needed to generate all relations in the semigroups
$S(m_1, m_2, m_3)$.
Then the semigroup $U_{c(\tree)}^L(\br)$ has relations
generated in degree bounded by $N$. 
\end{proposition}
This proposition, coupled with Theorem \ref{cuberel}
proves Theorem \ref{rel}. We recap the content of the
last two sections with the following theorem. 

\begin{theorem}
Let $(\tree, \br, L)$ be admissible.  Then
the ring $\C[U_{c(\tree)}^L(\br)]$ 
has a presentation
$$
\begin{CD}
0 @>>> I @>>> \C[X] @>>> \C[U_{c(\tree)}^L(\br)] @>>> 0\\
\end{CD}
$$
where $X$ is the set of degree $1$ elements of
$U_{c(\tree)}^L(\br)$, and $I$ is the ideal generated
by two types of binomials, 
$$[\omega_1]\circ \ldots \circ[\omega_n] - [\eta_1]\circ \ldots \circ[\eta_n].$$
\begin{enumerate}
\item Binomials where $n \leq 3$, $i_v^*(\omega_1) + \ldots + i_v^*(\omega_n) = i_v^*(\eta_1) + \ldots + i_v^*(\eta_n)$
is a balanced relation in $U_Y^L$ for some specific $v$, and 
$\{i_{v'}^*(\omega_1), \ldots, i_{v'}^*(\omega_n)\} = \{i_{v'}^*(\eta_1), \ldots, i_{v'}^*(\eta_n)\}$
for $v \neq v'$. 
\item Binomials where $n = 2$ and $i_v^*(\omega_1) + i_v^*(\omega_2) = i_v^*(\eta_1) + i_v^*(\eta_2)$
such that $\{i_v^*(\omega_1), i_v^*(\omega_2)\}$ is balanced for all $v \in I(c(\tree))$.
\end{enumerate}
This induces a presentation for $\C[S_{\tree}^L(\br)]$ by isomorphism.
\end{theorem}

\begin{corollary}
The same holds for $\C[S_{\tree}(\br)]$.
\end{corollary}

\begin{proof}
For each pair $(\tree, \br)$ 
it is easy to show that there is a number $N(\tree, \br)$, 
such that any weighting $\omega$ which satisfies
the triangle inequalities on $\tree$
and has $\omega(e_i) = \br_i$ must have
$\omega(e) \leq N(\tree, \br)$ for $e \in E(\tree).$
Because of this $S_{\tree}^L(\br) = S_{\tree}(\br)$
for $L$ sufficiently large.
\end{proof}

\section{Special Cases and Observations}

In this section we collect results on 
some special cases of $\C[S_{\tree}^L(\br)]$.
In particular we study some instances when cubic
relations are unnecessary, we give some examples
where the semigroup is not generated in degree $1$,
we analyze the case when $L$ is allowed to be odd, 
and we give instances where cubic relations are necessary.

\subsection{The Caterpillar Tree}
One consequence of the proof
of Theorem \ref{cuberel} is that
a semigroup $U_{c(\tree)}^{2L}(\br)$ which
omits or only partially admits the 
semigroup $S(0, 0, 0)$ or $S(L-1, L-1, 0)$
as an image of one of the morphisms $i_v^*$
manages to avoid degree $3$ relations entirely.  
The next proposition illustrates one
such example, the semigroups of weightings 
on the Caterpillar tree, pictured below.  

\begin{figure}[htbp]
\centering
\includegraphics[scale = 0.3]{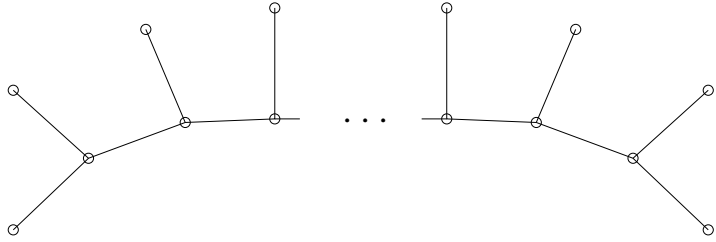}
\caption{The Caterpillar tree}
\label{fig:Cater}
\end{figure} 

\begin{proposition}
Let $\tree_0$ be the caterpillar tree, 
and let $\br(i)$ be even for all $i \in V(\tree_0)$, 
then $S_{\tree_0}^{2L}(\br)$  is generated
in degree 1, with relations generated by quadrics. 
\end{proposition}
 
\begin{proof}
We catalogue the
weights $i_v^*(\omega)$ which can appear in degree 1. 
For the sake of simplicity
we divide all weights by $2$.
Suppose $i_v(G)$ is an external edge, 
then $i_v^*(\omega)(E)$ and $i_v^*(\omega)(F)$
satisfy the following inequalities 

$$i_v^*(\omega)(E) \leq i_v^*(\omega)(F) + \frac{\br(i)}2$$
$$i_v^*(\omega)(F) \leq i_v^*(\omega)(E) + \frac{\br(i)}2$$
$$i_v^*(\omega)(E) + i_v^*(\omega)(F) + \frac{\br(i)}2 \leq 2L$$
where $i_v^*(\omega)(G) = \br(i)$. 
These conditions define a polytope in $\R^2$
with vertices $(L, L - \frac{\br(i)}2)$, 
$(L - \frac{\br(i)}2, L)$, $(\frac{\br(i)}2, 0)$
and $(0, \frac{\br(i)}2)$.  Pictured below
is the case $L = 9$, $\br(i) = 6$.
\begin{figure}[htbp]
\centering
\includegraphics[scale = 0.4]{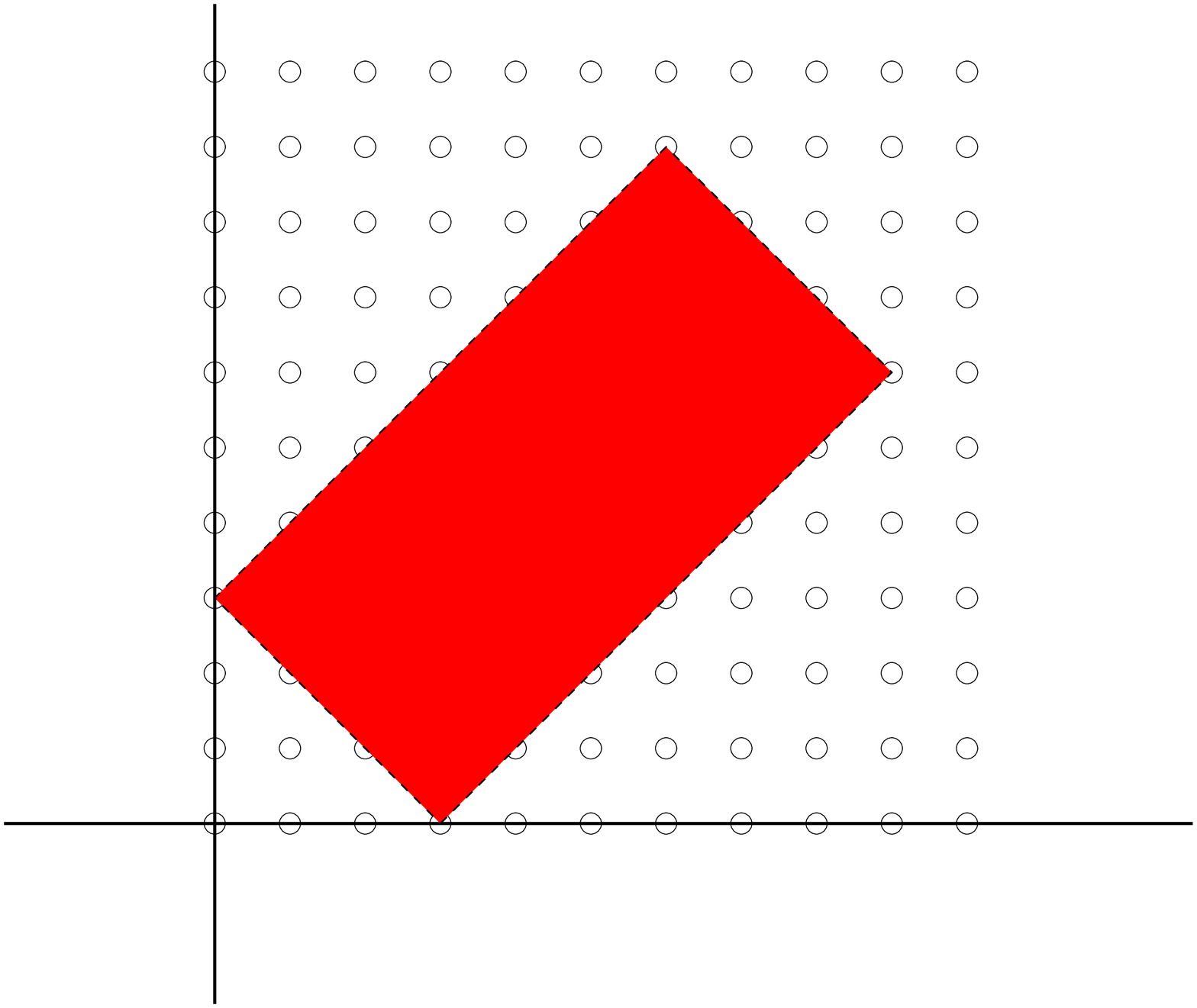}
\caption{The case $L = 9$, $\br(i) = 6$}
\label{fig:Graph}
\end{figure}
When two edges are external, the polytope
is an integral line segment. 
Note that the intersection of 
any lattice cube in $\R^2$ with the above
polytope is a simplex or a unit square. 
Both of these polytopes have at most
quadrics for relations in their semigroup
of lattice points.  Hence the argument
used to prove Theorem \ref{rel} shows
that $U_{c(\tree_0)}^{2L}(\br)$ needs 
only quadric relations.
\end{proof}

\begin{corollary}
If $L$ is even and greater than $2$, and
$\br$ is a vector of nonnegative even integers, 
the ring $R^G(L)_{\br}$ has a presentation
with defining ideal generated by quadrics. 
In particular, the second Veronese subring of
any $R^G(L)_{\br}$ has such a presentation
if $L > 1$. 
\end{corollary}

\subsection{Counterexamples to Degree $1$ generation}

Now we'll see how to generate
examples of $(\br, \tree, L)$ such that
$S_{\tree}^L(\br)$ is not
generated in degree $1$.  We will 
begin by defining a certain class
of paths in the tree $\tree$.

\begin{definition}
Let $\tree$ have an even
number of leaves.
Let $O(\tree)$ be the set of
paths in $\tree$ with the property that 
a weighting $\omega \in S_{\tree}$
which assigns all odd numbers to elements of $V(\tree)$, 
weights the edges of any member of $O(\tree)$ with 
an odd number under the parity condition. 
\end{definition}

Let us see that this is a well-defined set. 
It suffices to show that the parity of the 
members of $V(\tree)$ determines the parity
of every edge in $\tree$.  This follows
from induction on the number of edges in $\tree$. 
To see that members of $O(\tree)$ are paths
which never intersect, note that a lone
odd number can never appear in a trinode, 
nor can three odd numbers appear in a trinode. 
In particular, any pair of paired edges forms
a member of $O(\tree)$.   

\begin{figure}[htbp]
\centering
\includegraphics[scale = 0.3]{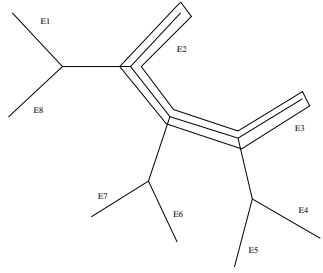}
\caption{E2 and E3 are lone leaves connected by an element of $O(\tree)$}
\label{fig:Lone Leaf}
\end{figure}

\begin{proposition}\label{nodeg1}
Let $(\br, \tree, L)$ be such
that the endpoints 
of each $\gamma \in O(\tree)$
are given the same parity, 
with some pair of endpoints
$(E, F)$ odd.  Then
if there is a degree $2$ weighting
which assigns $0$ to any edge
in $\gamma$, $S_{\tree}^{2L}(\br)$ is 
not generated in degree $1$. 
\end{proposition}

\begin{proof}
All degree $1$ elements
must assign odd numbers to
the edges on the path
joining $E$ and $F$.  No two 
odd numbers add to $0$. 
\end{proof}

\begin{corollary}\label{good}
The semigroup $S_{\tree}^{2L}(\vec{1})$ is generated
in degree $1$ if and only if $\tree$ has the
property that no leaf is lone and $L > 1$
\end{corollary}

\begin{proof}
The condition that the members of 
$\br$ sum to an even number forces
us to only consider trees $\tree$
with an even number of leaves. 
First we show that a tree with
lone leaves has a degree $2$ 
weighting satisfying the conditions
of proposition \ref{nodeg1}. 
Since $L > 1$, it suffices to note 
that for any tree $\tree$, and internal
edge $e \in \tree$, there is a weighting
that assignes the edge $e$ zero and every
other edge $2$.
If $\tree$ contains only paired leaves,
we can restrict to the tree $c(\tree)$ and consider
halved weightings without the parity condition. 
In this context, the weighting which assigns
every edge $1$ can be factored only if 
$L > 1$.  This finishes the $only$ $if$ portion of the statement. 
The $if$ portion of the statement is taken care of by Theorem \ref{gen}. 
\end{proof}

\begin{remark}
Trees with the property that no leaf is lone are called
Good Trees in \cite{HowardMillsonSnowdenVakil}, 
where they were introduced by Andrew Snowden for the purpose
of proving the analogue of Corollary \ref{good} for 
$S_{\tree}(\vec{1})$.
\end{remark}

\subsection{The Case when $L$ is odd}

When the level $L$ is odd, the polytope
$P_3(L)$ is no longer integral, however its
Minkowski square $P_3(2L)$ is integral, so clearly
there are elements of $P_3(2L)$ which cannot
be integrally factored, specifically the corners. 
This observation has a generalization.

\begin{definition}
Let $IP_3(L)$ be the convex hull of
the integral points of $P_3(L)$.
Let $\Omega$ be the set of elements in
the graded semigroup of lattice points of $P_3(L)$
such that $\frac{1}{deg(Q)}Q \in P_3(L) \setminus IP_3(L)$. 
\end{definition}

Let $(E, F, G) = Q \in P_3(L)$ be integral with $L$ odd, and
suppose $E$, $F$, or $G \geq \frac{L-1}2 + 1$.
Then, by the triangle inequalities we must have
$F + G \geq \frac{L-1}2 + 1$, so $E + F + G \geq L+1$, 
a contradiction.  This shows that $IP_3(L)$ is
contained in the intersection of $P_3(L)$ with the
halfspaces $E, F, G \leq \frac{L-1}2$, this identifies
$IP_3(L)$ as the convex hull of the set
$$\{(0,0,0), (\frac{L-1}2, \frac{L-1}2, 0),  (\frac{L-1}2, 0, \frac{L-1}2),(0, \frac{L-1}2, \frac{L-1}2),$$
$$(\frac{L-1}2, \frac{L-1}2, 1),  (\frac{L-1}2, 1, \frac{L-1}2),(1, \frac{L-1}2, \frac{L-1}2)\}.$$
The case $IP_3(5)$ is pictured below. 

\begin{figure}[htbp]
\centering
\includegraphics[scale = 0.4]{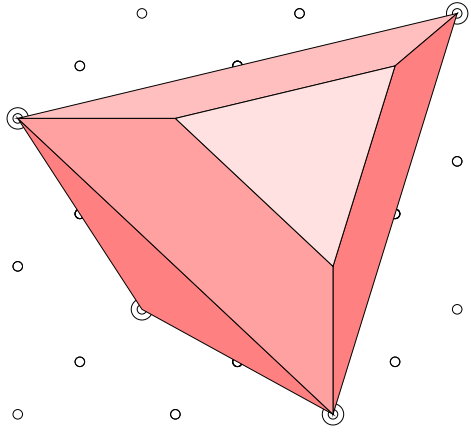}
\caption{The Polytope $IP_3(5)$}
\label{IPfigure}
\end{figure}

\begin{proposition}
Any $Q \in \Omega$ cannot be integrally factored. 
\end{proposition}

\begin{proof}
This follows from the observation
that if $Q = E_1 + \ldots + E_n$
then $\frac{1}{n}Q$ is in the convex 
hull of $\{E_1, \ldots, E_n\}$.
\end{proof}
A factorization
of any element $\omega$ such that $i_v^*(\omega) = Q$
gives a factorization of $Q$.  So any
$\omega \in U_{c(\tree)}^L(\br)$ with
a $i_v^*(\omega) \in \Omega$ is necessarily
an obstruction to generation in degree $1$, 
this also turns out to be a sufficient obstruction
criteria. 

\begin{theorem}
Let $\tree$ and $\br$ satisfy the same
conditions as admissibility, and let $L \neq 2$.  Then
$U_{c(\tree)}^L(\br)$ is generated in degree
$1$ if and only if 

$$i_v^*(\omega) \in U_Y^L \setminus \Omega$$
for all $v \in I(c(\tree))$, $\omega \in U_{c(\tree)}^L(\br)$.
In this case all relations
are generated by those of degree
at most $3$. 
\end{theorem}

\begin{proof}
We analyze $IP_3(L)$ in the same way we did $P_3(2L)$.
The reader can verify that the integral
points of $C(m_1, m_2, m_3) \cap P_3(L)$ are the
same as the integral points of $C(m_1, m_2, m_3) \cap IP_3(L)$.
The possibilities are represented
by slicing the cubes in Figure \ref{fig:Parity} along
the plane formed by the upper right or lower left
collection of three non-filled dots, depending on the cube,
and then restricting to the convex hull of the remaining integral
points. All cases are lattice
equivalent to one of the polytopes
listed in Figure \ref{fig:Cube}, after considering
two and one dimensional cases as facets of
neighboring three
dimensional polytopes. Since any element of $U_Y^L$
not in $\Omega$ is necessarily a lattice point of a Minkowski
sum of $IP_3(L)$, 
the theorem follows by the same arguments
used to prove Theorems \ref{gen} and \ref{rel}
\end{proof}

\subsection{Necessity of Degree $3$ Relations}

Now we show that there are
large classes of admissible $(\tree, \br, L)$ which
require degree $3$ relations. We will exhibit a degree $3$ 
weighting which has only two factorizations. The tree $\tree$ with weight 
$\omega_{\tree}$ is pictured below, it is an element
of $S_{\tree}(\vec{2})$. In all 
that follows all weightings are considered to have been halved.

\begin{figure}[htbp]
\centering
\includegraphics[scale = 0.4]{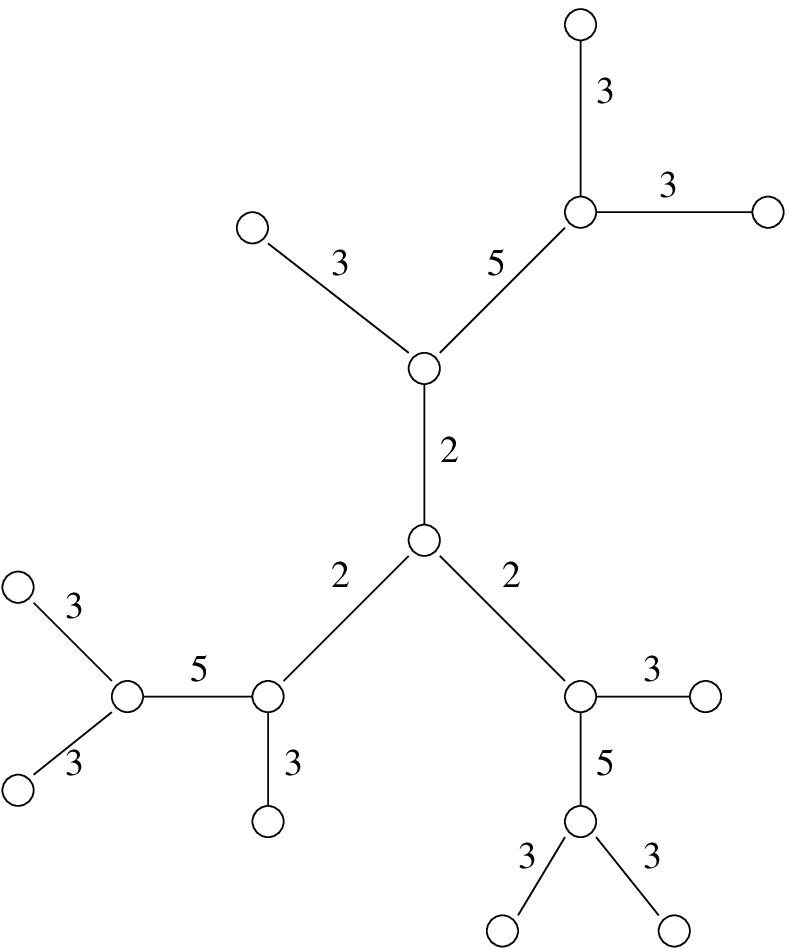}
\caption{$\omega_{\tree}$}
\label{}
\end{figure}

\begin{figure}[htbp]
\centering
\includegraphics[scale = 0.4]{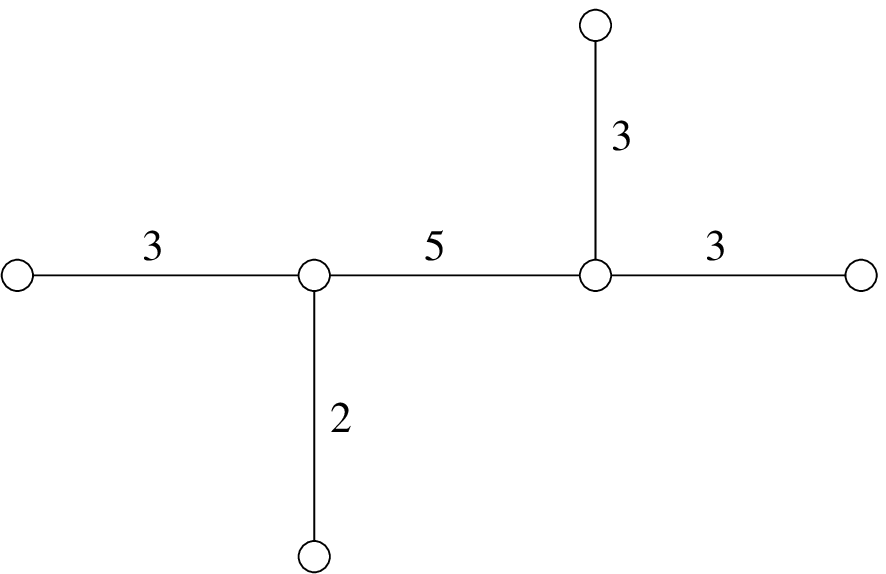}
\caption{ $\omega_{\tree'}$}
\label{NC1}
\end{figure}

\begin{figure}[htbp]
\centering
\includegraphics[scale = 0.4]{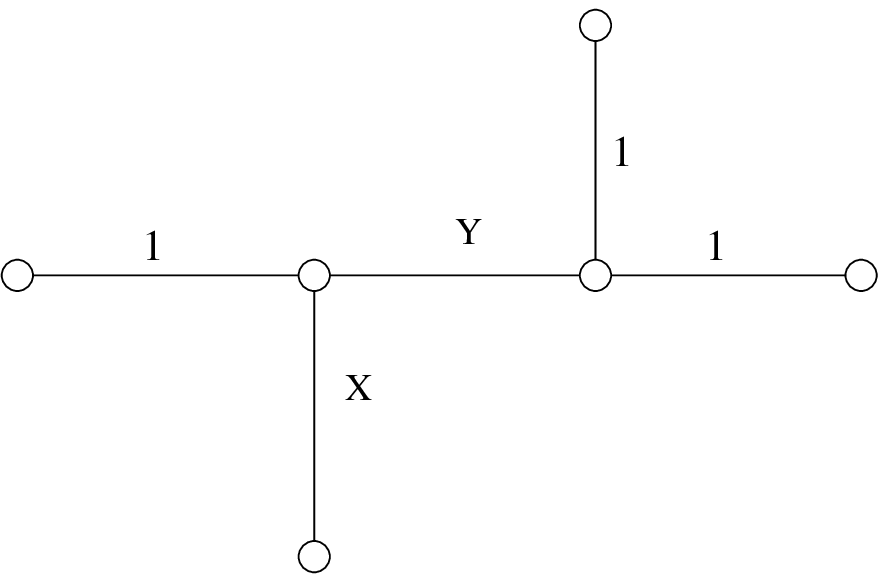}
\caption{}
\label{counter3}
\end{figure}

Notice that $\omega_{\tree}$ has 3-way symmetry 
about the central trinode, we will exploit this
by considering the tree $\tree'$ with 
restricted weighting $\omega_{\tree'}$ pictured in
Figure \ref{NC1}.  We find the 
weightings that serve as a degree 1 factors of
$\omega_{\tree'}$.  First of all, any degree 1 
weighting which divides $\omega_{\tree'}$ must be as in 
Figure \ref{counter3}.

\begin{figure}[htbp]
\centering
\includegraphics[scale = 0.4]{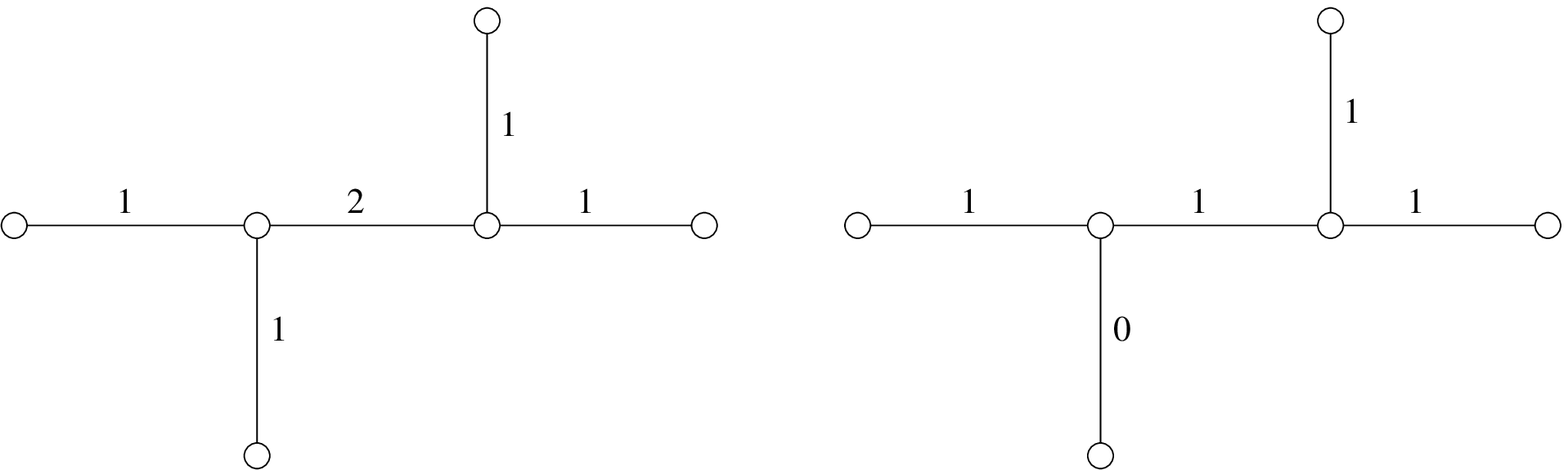}
\caption{}
\label{NC2}
\end{figure}

It suffices to find the possible values of $X$ and $Y$.  
Both must be $\leq 2$, which shows that $Y$ can be either
$2$ or $1$.  Now, by the triangle inequalities, any $X$
paired with $Y = 1$ must be $\geq 1$. Since two members
of a factorization must have $Y = 2$, $X$ must also
have a value $\leq 1$ on these factors. There are
exactly two possibilities determined by the value of $X$,
both are shown in Figure \ref{NC2}. Any factorization of $\omega_{\tree}$ is determined
by its values on the central trinode, and these values must be weights
composed entirely of $0$ and $1$.  There are exactly two such variations, 
making (of course) the Degenerated Segre Cubic.

We have not specified a level $L$ for this weighting, but the same 
argument applies for any level 
large enough to admit $\omega_{\tree}$ as a weighting in degree $3$.
For any tree $\tree^*$, edge $e^* \in tree^*$, and weight $\omega_{\tree^*}$ we can 
create a new weight on a larger tree by adding a vertex in the middle of $e^*$,
attaching a new leaf edge at that vertex, and weighting the both sides of
the split $e^*$ with $\omega_{\tree^*}(e^*)$, and the new edge with $0$. Using 
this procedure on any $(\tree^*, e^*, \omega_{\tree^*})$, and $(\tree, e, \omega_{\tree})$
for any edge $ e \in \tree $, can create a new weighted tree by identifying
the new $0$-weighted edges. An example of this procedure, which we call 
merging, is pictured below. In this way many examples of unremoveable
degree $3$ relations can be made. 

\begin{figure}[htbp]
\centering
\includegraphics[scale = 0.3]{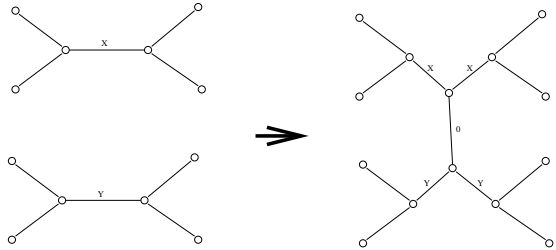}
\caption{Merging two tree weightings.}
\label{merge}
\end{figure}

\bigskip

Christopher Manon:
Department of Mathematics, 
University of Maryland, 
College Park, MD 20742, USA, 
manonc@math.umd.edu

\end{document}